\documentclass[12pt]{article}
\usepackage[default]{jasa_harvard}    
\usepackage{JASA_manu}
\usepackage{setspace}
\usepackage{url} 
\usepackage{amsmath,graphicx,amssymb,amsfonts,subfigure,mathtools}
\usepackage{enumerate}
\usepackage{bm}
\usepackage{color, colortbl}
    \usepackage{algorithmicx}
    \usepackage[ruled]{algorithm}
    \usepackage{algpseudocode}
    \algnewcommand\algorithmicinput{Input:}
	\algnewcommand\INPUT{\item[\algorithmicinput]}
	\algnewcommand\algorithmicoutput{Output:}
	\algnewcommand\OUTPUT{\item[\algorithmicoutput]}

\usepackage{multirow}
\usepackage{mathtools}
\usepackage{array}
\usepackage[table]{xcolor}

\usepackage{float}
\floatstyle{plain}
\restylefloat{table}

\def\AIC{\textsc{aic}}
\def\BC{\textsc{bc}}
\def\HQ{\textsc{hq}}
\def\BIC{\textsc{bic}}

\def\PI{\textsc{pi}_N}
\def\Beta{\mathcal{B}}

\def\T{{ \mathrm{\scriptscriptstyle T} }}
\def\P{{ \mathrm{pr} }}
\def\v{{\varepsilon}}
\def\LM{{L_{\textrm{max}}^{(N)}}}
\def\lm{{L_{\textrm{max}}}}
\def\M{M_N}

\def\Laic{\hat{L}_{\textsc{aic}}}
\def\Lbic{\hat{L}_{\textsc{bic}}}
\def\Lbc{\hat{L}_{\textsc{bc}}}
\def\Lopt{L_{0}^{(N)}}

\DeclareMathOperator*{\argmin}{arg\,min}
\newcommand{\rom}[1]{\uppercase\expandafter{\romannumeral #1\relax}}
\DeclarePairedDelimiterX{\norm}[1]{\lVert}{\rVert}{#1}
\DeclarePairedDelimiterX{\abs}[1]{\lvert}{\rvert}{#1}

\newtheorem{remark}{Remark}
\newtheorem{theorem}{Theorem}

\newtheorem{assumption}{Assumption}
\newtheorem{proposition}{Proposition}

\makeatletter
\renewcommand{\@makefnmark}{\hbox{\textsuperscript{\tiny{\@thefnmark}}}}
\makeatother

\begin{document}

\title{{\Large Bridging AIC and BIC: a new criterion for autoregression}}

\author{
Jie Ding,  
Vahid Tarokh \\
School of Engineering and Applied Sciences, Harvard University \\ 
Yuhong Yang \\
School of Statistics, University of Minnesota \\ 
}

\date{\vspace{-3ex}}

\maketitle


\begin{abstract}
We introduce a new criterion  to determine the order of an autoregressive model fitted to time series data.
It has the benefits of the two well-known model selection techniques, the Akaike information criterion and the Bayesian information criterion.
When the data is generated from a finite order autoregression, 
the Bayesian information criterion is known to be consistent, and so is the new criterion.
When the true order is infinity or suitably high with respect to the sample size, the Akaike information criterion is known to be efficient in the sense that its prediction performance is asymptotically equivalent to the best offered by the candidate models; in this case, the new criterion behaves in a similar manner.
Different from the two classical criteria, the proposed criterion adaptively achieves either consistency or efficiency depending on the underlying true model.
In practice where the observed time series is given without any prior information about the model specification, the proposed order selection criterion is more flexible and robust compared with classical approaches.  
Numerical results are presented demonstrating the adaptivity of the proposed technique when applied to various datasets.
\end{abstract}



\section{Introduction}
\label{sec:Intro}


In a practical situation of the autoregressive model fitting,  the order of the model is generally unknown.
There have been many  order selection methods proposed, following different philosophies.
Anderson's multiple decision procedure \cite{anderson1962determination} 
sequentially tests when the partial autocorrelations of the time series become zero.
The final prediction error criterion proposed by \citeasnoun{akaike1969fitting} aims to minimize the one-step prediction error when the estimates are applied to another independently generated dataset.
 \citeasnoun{bhansali1977some} generalized the final prediction error criterion  by replacing  $2$ with a parameter $\alpha$ in its formula, and 
proved that the asymptotic probability of choosing the correct order increases as  $\alpha$ increases.
The well-known Akaike Information Criterion, AIC \cite{akaike1998information}, was derived by minimizing the Kullback-Leibler divergence between the true distribution and the estimate of a candidate model.
Some variants of AIC, for example the modified Akaike information criterion that replaces the constant $2$ by a different positive number, have also been considered \cite{broersen2000finite}.
Nevertheless, \citeasnoun{akaike1979bayesian} argued in a Bayesian setting that the original AIC is more reasonable than its variants in a practical situation.
\citeasnoun{hurvich1989regression}  proposed the corrected AIC for the case where the sample size is small.
Another popular method is the Bayesian information criterion, BIC, proposed by \cite{schwarz1978estimating} that aims at selecting a model that maximizes the posterior model probability.
\citeasnoun{hannan1979determination}  proposed a criterion, HQ, that replaces the $\log N$ term in BIC by $c \log \log N (c>1)$, where $N$ is the sample size, and they showed that it is the smallest penalty term that guarantees strong consistency of the selected  order.
The focused information criterion is another remarkable approach that takes into account the specific purpose of the statistical analysis, by estimating the risk quantity of interest for each candidate model \cite{claeskens2003focused,claeskens2007prediction}.
%
%
%
Other methods for autoregressive order selection include
the criterion autoregressive transfer function method \cite{parzen1974some},
the predictive least-squares principle \cite{rissanen1986predictive,hemerly1989strong},
the combined information criterion \cite{broersen2000finite};
see  \citeasnoun{de1985methods} and \citeasnoun{shao1997asymptotic} for more references. 
Despite the rich literature on autoregressive models, the most common order selection criteria are AIC and BIC.


In this paper, the specified model class for fitting is the set of autoregressions with orders $L=1,\ldots,\lm$ for some prescribed natural number $\lm$.  
In relation to the true data generating process, the model class is referred to as \textit{well-specified} (or parametric) if the data is generated from a finite order autoregression and the true order is no larger than $\lm$, and \textit{mis-specified} (or nonparametric) if otherwise.
It is well known that BIC is consistent in order selection in the well-specified setting.
In other words, the probability of choosing the true order tends to one as the sample size tends to infinity.
The Akaike information criterion is not consistent and has a fixed overfitting probability when the sample size tends to infinity \cite{shibata1976selection}.
However, AIC is shown to be efficient in the mis-specified setting, while BIC is not \cite{shibata1980asymptotically}. 
Here we call an order selection procedure (asymptotically) efficient if its prediction performance (in terms of the squared difference between the prediction and its target conditional mean) is asymptotically equivalent to the best offered by the candidate autoregressive models. A rigorous definition of efficiency is given in Section~\ref{sec:efficiency}. 
In other words, AIC typically produces less modeling error than BIC when the data is not generated from a finite order autoregressive process. Furthermore, asymptotic efficiency of AIC for order selection in terms of the same-realization predictions for infinite order autoregressive or integrated autoregressive processes has also been well established \cite{ing2005orderselection,ing2012integratedseries}.

In real applications, one usually does not know whether the model class is well-specified. The task of adaptively achieving the better performance of AIC and BIC is theoretically intriguing and practically useful.
There have been several efforts towards this direction. 
\citeasnoun{yang2005can} considered the possibility of sharing the strengths of AIC and BIC in the regression context.
It has been shown under mild assumptions that any consistent model selection criterion behaves suboptimally for estimating the regression function in terms of the \textit{minimax} rate of convergence.
In other words, the conflict between AIC and BIC in terms of achieving model selection consistency and minimax-rate optimality in estimating the regression function cannot be resolved.  
But this does not indicate that there exists no criterion achieving the pointwise asymptotic efficiency in both well-specified and mis-specified scenarios, because the minimaxity (uniformity over the linear coefficients) is intrinsically different from the (pointwise) efficiency.
In the remarkable work by \citeasnoun{ing2007accumulated}, a hybrid selection procedure combining AIC and a BIC-like criterion was proposed. Loosely speaking, if a BIC-like criterion selects the same model at  sample sizes $N^{\ell}$ ($0<\ell<1$) and $N$, then with high probability (for large $N$) the model class is well-specified and the true model has been converged to, and thus a BIC-like criterion is used; otherwise AIC is used. Under some conditions, the hybrid criterion was proved to achieve the pointwise asymptotic efficiency in both well-specified and mis-specified scenarios.
In estimating regression functions with independent observations, \citeasnoun{yang2007prediction} proposed a similar approach to adaptively achieve asymptotic efficiency for both parametric and nonparametric  situations, by examining whether BIC  selects the same model again and again at different sample sizes (instead of only two sample sizes used by \citeasnoun{ing2007accumulated}).
\citeasnoun{liu2011parametric} proposed a method to adaptively choose between AIC and BIC based on a measure called parametricness index.
In the context of sequential Bayesian model averaging, \citeasnoun{erven2012catching} and \citeasnoun{pas2015catchingup} used a switching distribution to encourage early switch to a better model and offered interesting theoretical understanding on its simultaneous properties.
Cross-validation has also been proposed as a general solution to choosing between AIC and BIC. It was shown by \citeasnoun{zhang2015cross} that, with a suitably chosen data splitting ratio, the composite criterion asymptotically behaves like the better one of AIC and BIC for both the AIC and BIC territories. 

In this paper, we introduce a new model selection criterion which is referred to as the bridge criterion (BC) for autoregressive models.
The bridge criterion is able to address the following two issues: First, given a realistic time series data, an analyst is usually unaware of whether the model class is well-specified or not; Second, even if the model class is known to be correct, the order (dimension) is not known, so that any prescribed finite candidate set suffers the risk of missing the true model.
We show that BC achieves both consistency when the model class is well-specified and asymptotic efficiency when the model class is mis-specified under some sensible conditions. 
Recall that the penalty terms of AIC and BIC are proportional to $L$  for autoregressive model of order $L$. 
In contrast, a key element of BC is the expression  $1+2^{-1}+\cdots+L^{-1}$ employed in its penalty term. As we shall see, it is the harmonic number that ``bridges'' the features of AIC and BIC.
Another key element is to let $\lm$ grow with sample size.
We emphasize that for the well-specified case, once the true order is selected with probability close to one, the resulting predictive performance is also asymptotically optimal/efficient. From this angle, the criterion achieves the asymptotic efficiency for both the well-specified and the mis-specified cases.

The outline of this paper is given below.
In Section~\ref{sec:Formulation}, we formulate the problem considered in this paper and briefly introduce the background and how the new criterion was heuristically derived.
In Section~\ref{sec:Lag_Selection}, we propose the bridge criterion and give an intuitive interpretation of it. 
We establish the consistency and the asymptotic efficiency property in Section~\ref{sec:theory}.
Numerical results are given in Section~\ref{sec:Num_Results} comparing the performance of our approach and other techniques.
In Section~\ref{sec:modifiedBC}, we propose a two-step strategy to adaptive choose the candidate size $\lm$, in order to further relax the conditions required by the theorems established in previous sections.  To that purpose, we also extend the expression of the bridge criterion. 
Finally, we make some discussions in Section~\ref{sec:Conclusion}.
%
 


\section{Background} \label{sec:Formulation} 

\subsection{Notation} \label{subsec:Notations}

We use $o_p(1)$ and $O_p(1)$ to denote any random variable that converges in probability to zero, and that is stochastically bounded, respectively.
We write $h_N=\Theta(g_N)$ if $c <  h_N/g_N < 1/c$ for some positive constant $c$ for all sufficiently large $N$, and $h_N=O(g_N)$ if $|h_N| < c g_N$ for some positive constant $c$ for all sufficiently large $N$.
If $\lim_{N \rightarrow \infty} f_N / g_N = 0 $, we write $f=o_N(g)$, or for brevity, $f=o(g)$.
Let $\lfloor x \rfloor$ denote the largest integer  less than or equal to $x$.
Let $\mathcal{N} (\mu, \sigma^2), \Beta (a,b), \chi_k^2$ respectively denote the normal distribution with density function
$f(x) = \exp \{ -(x-\mu)^2 / (2 \sigma^2) \} / (\sqrt{2 \pi } \sigma) $,
the Beta distribution with density function $f(x) = x^{a-1} (1-x)^{b-1} / B(a,b)$, where $B(\cdot,\cdot)$ is the beta function,
and the chi-square distribution with $k$ degrees of freedom.

\subsection{Problem formulation}

Given observations $\{x_n: n=1,\ldots,N_0\}$, we consider the following autoregressive model of order $L \ (L \in \mathbb{N})$ 
    \begin{align} \label{AR}
         x_n + \psi_{L,1} x_{n-1} + \cdots + \psi_{L,L} x_{n-L}  =\epsilon_n,
    \end{align}
where 
$\psi_{L, \ell} \in \mathbb{R}$ ($\ell=1,\ldots,L$), $\psi_{L, L} \neq 0$, the roots of the polynomial $z^{L} + \sum_{\ell=1}^L \psi_{L, \ell} z^{L-\ell}$ have modulus less than 1, and $\varepsilon_n$'s are independent and identically distributed according to $\mathcal{N}(0,\sigma^2)$.
%
%
%
The autoregressive model is referred to as AR$(L)$ model, and $ [ \psi_{L,1}, \ldots, \psi_{L,L}]^{\T}$ is referred to as the stable autoregressive filter $\Psi_{L} $.
Let $L_0$ denote the true order, which is considered to be finite for now.
In other words, the data is generated in the way described by (\ref{AR}) with $L=L_0$.
When $L_0$ is unknown, we assume that $\{1,\ldots, \lm\}$ is the candidate set of orders.
Let $N = N_0 - \lm$. 
The sample autocovariance vector and matrix  are respectively $ \hat{ \gamma}_{L}= [\hat{\gamma}_{1,0}, \ldots, \hat{\gamma}_{L,0}]^{\T} , \hat{ \Gamma}_L = [\hat{\gamma}_{i,j}]_{i,j=1}^L $.
where 
$
\hat{\gamma}_{i,j} = \frac{1}{N} \sum_{n=\lm+1}^{N_0} x_{n-i}  x_{n-j} \quad  (0 \leq i,j \leq \lm).
$
The filter of the autoregressive model of order $L$ can be estimated by
\begin{align} \label{LS}
\hat{ \Psi}_L = - \hat{ \Gamma}_L^{-1} \hat{ \gamma}_{L},
\end{align}
which yields consistent estimates \cite[Appendix 7.5]{box2011time}.
The one-step prediction error  is
%
 $       \hat{e}_L=    \sum_{n=\lm +1}^{N_0}  (x_n + \hat{\psi}_{L,1} x_{n-1} + \cdots + \hat{\psi}_{L,L} x_{n-L} )^2 / N.$
%
For convenience, we define $\hat{e}_0 = \hat{\gamma}_{0,0}$.
%
%
The error of the AR$(L)$ model can be calculated by
\begin{align}
    \hat{e}_L &=  \hat{e}_0 - \hat{ \gamma}_{L}^{\T} \hat{ \Gamma}_L^{-1} \hat{ \gamma}_{L}. \label{relation2}
\end{align}
Let $\gamma_{i-j}= E \{ x_{n-i}  x_{n-j} \}$ ($i,j \in \mathbb{Z}$) be the autocovariances and $ \Psi_L = [ \psi_{L,1}, \ldots, \psi_{L,L}]^{\T}$ be the best linear predictor of order  $L$.
In other words, $\Psi_{L}$ ($L \geq 1$) is the minimum of 
\begin{align}
  e_L = &  \min_{\psi_{L,1}^{*},\ldots,\psi_{L,L}^{*} \in \mathbb{R}} E  \bigl\{ (x_n + \psi_{L,1}^{*} x_{n-1} + \cdots + \psi_{L,L}^{*} x_{n-L} )^2 \bigr\},\label{var}
\end{align}
where the expectation is taken with respect to the stationary process $\{X_n\}$. In addition, we define $e_0 = \gamma_{0}$.
The values of $\Psi_{L} $ and $e_L$ can be calculated from a set of equations similar to (\ref{LS})--(\ref{relation2}), by removing the hats ($\wedge$) from all parameters.


Given an observed time series, the problem is how to identify the unknown order of the autoregressive model fitted to the data.
The Akaike information criterion and the Bayesian information criterion for autoregressive order selection is to select $\hat{L}$ ($1 \leq \hat{L} \leq \lm$) that respectively minimizes the  quantities
%
%
$ \AIC (N,L) = \log \hat{e}_L + 2 L/N,  \,
        \BIC (N,L) = \log \hat{e}_L +  L \log(N) / N $.
In the following two subsections, we introduce the motivation and perspective that naturally led to the bridge criterion. The formal expression of BC and its performance in asymptotic regions are established in Sections~\ref{sec:Lag_Selection} and ~\ref{sec:theory}.

\subsection{Motivation}

Distinct from AIC or BIC, the new criterion was initially derived from some perspectives unique to autoregressions.
Briefly speaking, it was initially motivated by postulating that nature randomly draws the coefficients of true autoregressions from a non-informative uniform distribution 
and by fixing the
type \rom{1} error 
in a sequence of hypothesis tests on the order. 
Suppose that we generate a time series to simulate an AR$(L_0)$ process using (\ref{AR}).
Clearly,
$\hat{e}_1 \geq \cdots \geq \hat{e}_{L_0-1} \geq \hat{e}_{L_0} \geq \hat{e}_{L_0+1} \geq \cdots \geq \hat{e}_{\lm}.$ 
Because of  (\ref{relation2}) and the consistency of $\hat{ \psi}_L$,  generally $\hat{e}_L$ is large for $L < L_0$ and is  much smaller for $L \geq L_0$.
If we plot $\hat{e}_L$ against  $L$ for $L=1,\ldots,\lm$, the curve is usually decreasing for $L < L_0$ and becomes almost flat for $L>L_0$. Intuitively, the order $\hat{L}$ may be selected such that $\hat{e}_L/ \hat{e}_{L-1}$ becomes ``less significant'' than its predecessors for $L > \hat{L} $.
%
%
%
We define the empirical and theoretical gain of goodness of fit using AR$(L)$ over AR$(L-1)$, respectively, as
\begin{align}
 \hat{g}_{L} &= \log \biggl(\frac{ \hat{e}_{L-1} }{ \hat{e}_{L}} \biggr) ,
 \quad 
 g_{L} = \log \biggl( \frac{e_{L-1} }{ e_{L} } \biggr) .
 \label{gain_2}
 \end{align}

Suppose that the data is generated by a stable filter $ \Psi_{L_0}$ of order  $L_0$.
For any positive integer $L $ that is greater than $L_0$ and does not depend on $N$, it was shown by \citeasnoun[Theorem 5.6.2 and 5.6.3]{anderson2011statistical} that $\sqrt{N}  [\hat{\psi}_{L_0+1,L_0+1}, \ldots,\hat{\psi}_{L ,L }]^{\T}$ has a limiting joint-normal distribution $\mathcal{N}( 0,  I)$ as $N$ tends to infinity, where $ I$ denotes the identity matrix.
In addition, the random variables $N \hat{g}_{L } \ (L =L_0+1,\ldots,\lm)$ are asymptotically independent and distributed according to $\chi_1^2$, where $\lm > L_0$ is a constant that does not depend on~$N$ \cite{shibata1976selection}. 
%
Next, we revisit AIC and BIC by associating them with a sequence of hypothesis tests. The purpose of the argument below is to motivate our new criterion. 

Test: We choose a fixed number $0<q<1$ as the significance level (or the type \rom{1} error), and thresholds $s$ such that
$q = \P ( W  > s)$,  
where $W \sim \chi_1^2$.
Consider the hypothesis test
\begin{align} \label{hypothesis1}
&H_0: L_0= 
		L -1 \quad
H_1:  L_0\geq  
		L .
\end{align}
If $N\hat{g}_L > s$ (or equivalently $s/N-\hat{g}_L < 0$), we reject $H_0$  and  replace $L-1$ by $L$, for $L = 2,3,\ldots$  until $L = \lm$ or $H_0$ is not rejected. 
One limitation of this hypothesis test technique is that it may produce extreme values \cite{akaike1970statistical}. 
A straightforward alternative solution would be to select the $L$ such that the aggregation of $s/N-\hat{g}_1,\ldots,s/N-\hat{g}_L $ is minimized, i.e., to select the global minimum:
\begin{align}\label{forward}
\hat{L} = \argmin_{1 \leq L \leq \lm} \sum\limits_{k=1}^L \bigl(\frac{s}{N}-\hat{g}_k \bigr) =  \log \hat{e}_L+ \frac{sL}{N}   - \log \hat{e}_0,
\end{align}
the objective function of which can be regarded as the goodness of fit $\hat{e}_L$ plus the penalty of the model complexity. 
The penalty term is a sum of thresholds $s$ and  $-\log \hat{e}_0$. The term $-\log \hat{e}_0$ does not depend on $L$, so it has no effect on the produced result and is negligible. 
The Akaike information criterion has a penalty term $2L/N $, it therefore corresponds to the above hypothesis tests with $q=$0.1573 . 
The Bayesian information criterion has a penalty term $L \log(N)/N$. It corresponds to the hypothesis tests with varying $q$.
As an illustration, the significance levels $q$ of BIC under different sample sizes are tabulated in Table~\ref{table:BIC}.

\begin{table}[htb]
\caption{\label{table:BIC} Significance level $q$ of the Bayesian information criterion at different sample sizes}
\centering
\fbox{
\begin{tabular}{c c   c  c c  c   }
 $N$  &100 &500 &1000 &2000 &10000 \\
  \hline
$q$  &0.0319  &0.0127 &0.0086 &0.0058 &0.0024 \\
\end{tabular}}
\end{table}

To motivate our new criterion, suppose that nature generates the data from an AR$(L_0)$ process, which is in turn randomly generated from the uniform distribution $\mathcal{U}_{L_0}$. 
Here, $\mathcal{U}_{L_0}$ is defined over the space of all the stable AR filters whose roots have modulus no larger than $r$ ($0 < r \leq 1)$: 
\begin{align} \label{space}
S_L(r) &=\biggl\{ \Psi_{L} : z^L + \sum_{\ell=1}^{L}\psi_{L, \ell} z^{L-\ell}
= \prod\limits_{\ell=1}^L (z-a_{\ell}) , \, \psi_{L,\ell} \in \mathbb{R}, |a_{\ell}| \leq r, \, \ell=1,\ldots,L \biggr\}.
\end{align}
Under this data generating procedure, $g_L$ is a random variable with distribution described by the following theorem. For the sake of continuity, we postpone a detailed discussion on $\mathcal{U}_{L_0}$ to the Supplementary Material.   

\begin{theorem} \label{thm:representation}
Suppose that $ \Psi_{L_0} $ is uniformly distributed in $S_L(1)$.
Then, $\psi_{1,1}, \ldots, \psi_{L_0,L_0}$ are independently distributed according to $(\psi_{L,L}+1)/2 \sim \Beta (\lfloor L/2+1 \rfloor,\lfloor (L+1)/2 \rfloor) \ (L=1,\ldots,L_0)$.
%
Furthermore, $L \psi_{L,L}^2 $  and $L  g_{L} $ converge in distribution to $\chi_1^2$ as $L$ tends to infinity.
\end{theorem}

Similarly, we postulate hypothesis tests in the opposite direction (for a given $\lm$):
\begin{align} \label{hypothesis2}
&H_0: 
L_0 = L \quad
H_1: 
L_0 \leq L-1 .
\end{align}
Under the null hypothesis, $g_L \neq 0$ almost surely, and we approximate the distribution of $\hat{g}_L$ by that of $g_L$.
We choose a fixed  number $0<p<1$ as the significance level,  and the associated thresholds $h_L$ at order $L$ such that $p = \P (g_L < h_L)$,
or equivalently
\begin{align} \label{conf_underfitting}
h_L = F^{-1}_{g_L} (p)
\end{align}
where $F^{-1}_{g_L}(\cdot)$ denotes the inverse function of the cumulative distribution function of $g_L$.
If $\hat{g}_L < h_L$ (or equivalently $\hat{g}_L-h_L<0$), we reject $H_0$ and replace $L$ by $L-1$, for $L = \lm, \lm - 1 , \ldots$ until $L=2$ or $H_0$ is not rejected.
Likewise, the $L$ that minimizes the following objective function can be chosen as the optimal order
\begin{align}
\hat{L} &= \argmin_{1\leq L \leq \lm} \sum\limits_{k=L+1}^{\lm} (\hat{g}_k - h_k)
= \log \hat{e}_L + \sum_{k=1}^{L} h_k + c \label{backward}
\end{align}
where $c=-(\log \hat{e}_{\lm} + \sum_{k=1}^{\lm} h_k  )$ does not depend on $L$.
The next subsection introduces the proposed criterion motivated by (\ref{backward}).


\subsection{Proposed order selection criterion}

From now on, we allow the largest candidate order to grow with the sample size $N$, and use notation $\LM$ instead of $\lm$ to emphasize this dependency. Define $N = N_0 - \LM$.
Building on the idea of  (\ref{backward}), 
we adopt the penalty term
$
 \sum_{k=1}^L h_k(p)
$
where $h_k(p)$ is defined in (\ref{conf_underfitting}), and $p$ is further determined by
\begin{align} \label{july_9_1}
h_{\LM}(p) = \frac{2}{N}.
\end{align}
%
Theorem~\ref{thm:representation} implies that $h_k(p) \approx F_{\chi_1^2}^{-1}(p)/k$ for large $k$, where $F_{\chi_1^2}^{-1}(\cdot)$ denotes the inverse function of the cumulative distribution function of $\chi_1^2$. From (\ref{july_9_1}) we have
$F_{\chi_1^2}^{-1}(p) \approx  2 \LM /N$, 
and thus
$ h_k(p)  \approx 2 \LM / (N k)$. 
%
We therefore propose the following bridge criterion:
select the $L \in \{1,\ldots,\LM\}$ that minimizes $\log \hat{e}_L + (2 \LM /N) \sum_{k=1}^L 1/k. $

We have seen that given a fixed type \rom{1} error, the threshold for hypothesis test (\ref{hypothesis1}) is a constant, while the threshold for (\ref{hypothesis2}) decreases in $L$ leading to the $1/k$ term. Intuitively speaking, the uniform distribution on $S_L(r)$ concentrates more around the boundary of the space, and the loss of underfitting, $e_{L-1}/e_{L}=1/(1-\psi_{L,L}^2)$, becomes more negligible, as $L$ increases.
To some extent, this observation suggests an interesting idea that the penalization for different models is not necessarily linear in model dimension; one may start with a BIC-type heavy penalty, but alleviate it more and more to an AIC-type light penalty as the candidate model is larger, offering the possibility of changing/reinforcing one's belief in the model specification.

\section{Bridge criterion} \label{sec:Lag_Selection}

\begin{figure} [t]
        \vspace*{-0.0 in}
	\begin{center}
    \hspace{-0.0 in}\includegraphics[width=5.0in]{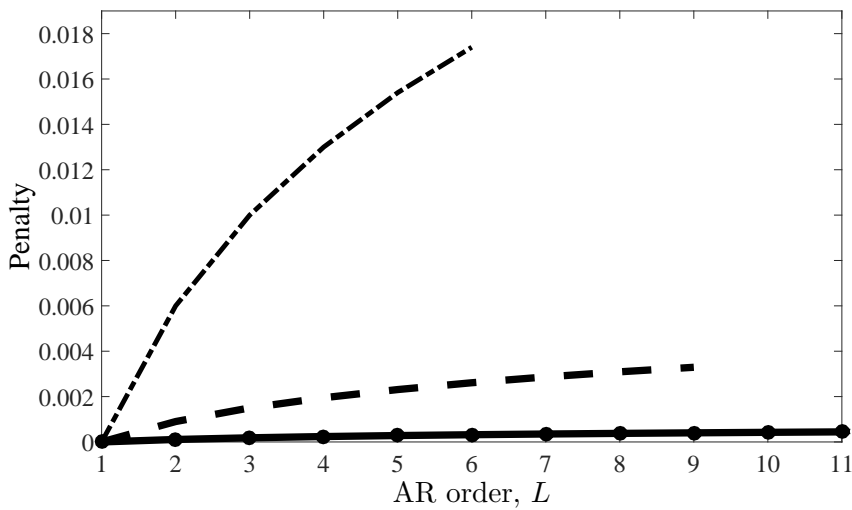}  
    \end{center}
    \vspace*{-0.2in}
    \caption{A graph showing the penalty term for sample size $10^3$ (dot-dash), $10^4$ (dashes), and $10^5$ (solid).}
    \label{fig_TS_samplesize}
    \vspace*{-0.21in}
    \end{figure}

Recall that 
the estimated order $\hat{L}$ by bridge criterion is
\begin{align}  \label{TS_criterion}
\hat{L} = \argmin_{1 \leq L \leq \LM} \BC(n, L)  = \log \hat{e}_L + \frac{2 \LM }{N} \sum_{k=1}^L \frac{1}{k}
\end{align}
where  $\LM$ is the largest candidate order.
$\LM$  must be selected such that $\lim_{N \rightarrow \infty}\LM = \infty$,
and its rate of growth will be studied in Section~\ref{sec:theory}.
It is well known that
 $\sum_{k=1}^L 1/k = \log L + c_E +o_L(1) $
for large $L$, where $c_E$ is the Euler-Mascheroni constant.
Fig.~\ref{fig_TS_samplesize} illustrates the penalty curves for different $N$ and  $\LM = \lfloor \log N \rfloor$.
Without loss of generality, we can shift the curves to be at the same position at $L=1$.


\begin{figure} [t]
    \vspace*{-0.0 in}
    \begin{center}
    \hspace{-0.0 in}\includegraphics[width=5.0in]{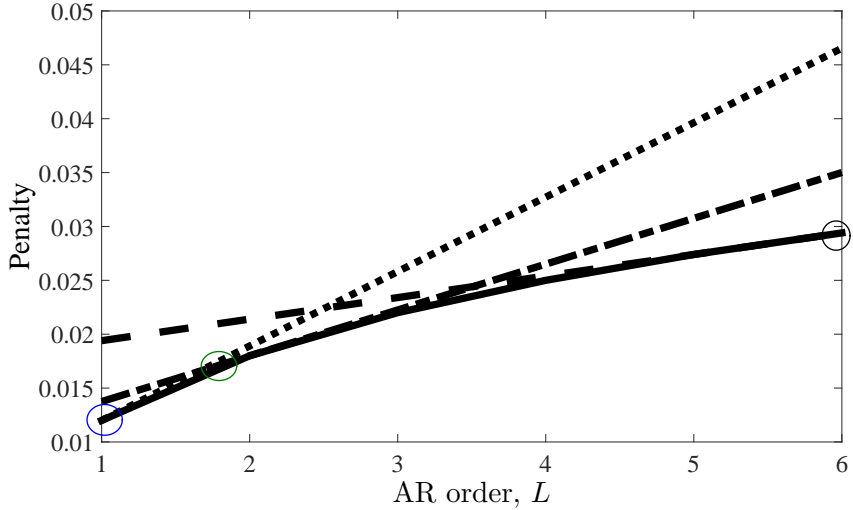}
    \end{center}
    \vspace*{-0.1 in}
    \caption{A graph showing the penalty curves of the bridge criterion (solid) together with the Akaike information criterion (dashes), the Hannan and Quinn criterion (dot-dash),  the Bayesian information criterion (small dashes), and the tangent points (circled) for $N=1000$ .  } 
    \label{fig_TS_comparison}
    \vspace*{-0.1 in}
    \end{figure}

Fig.~\ref{fig_TS_comparison} illustrates the penalty curves for the bridge criterion,  the Akaike information criterion, the Bayesian information criterion, and the Hannan and Quinn criterion, respectively denoted by
\begin{align*}
&J_{\BC}(L)=\frac{2 \LM}{N} \sum_{k=1}^L \frac{1}{k},\, 
J_{\AIC}(L)=\frac{2}{N}L , \,
J_{\BIC}(L)=\frac{\log(N) }{N} L,  \,
J_{\HQ}(L)=\frac{c\log \log(N)}{ N} L \, 
\end{align*}
where $c$ is chosen to be 1.1, $L=1,\ldots, \LM = \lfloor \log N \rfloor$, and $N=1000$.
Any of the above penalty curves can be written in the form of $\sum_{k=1}^L t_k$, and only the slopes $t_k \ (k=1,\ldots,\lm)$ matter to the performance of order selection.
For example, suppose that $L_2$ is selected instead of $L_1$ ($L_2 > L_1$) by some criterion. This implies that the gain of goodness of fit $\hat{e}_{L_1} - \hat{e}_{L_2}$ is greater than the sum of slopes $\sum_{k=L_1+1}^{L_2} t_k$.
Thus, we have shifted the curves of the latter three criteria to be tangent to the log-like curve of the bridge criterion in order to highlight their differences and connections.
Here, two curves are referred to as tangent to each other if they intersect at one and only one point, the tangent point.
The tangent points (marked by circles) of $J_{\AIC}$, $J_{\HQ}$ and $J_{\BIC}$ are respectively 6, 2 and 1.
Take the curve $J_{\HQ}$ as an example. The meaning of the tangent point is that BC penalizes more than HQ  for $k \leq 2$ and otherwise for $k > 2$.

Given a sample size $N$, the tangent point between $J_{\BC}$ and $J_{\HQ}$ curves is at $T_{\BC:\HQ}=2\LM / (c \log \log N) $.
As an example, we choose $\LM = \lfloor \log N \rfloor$.
If the true order $L_0$ is finite, $T_{\BC:\HQ}$  will  be larger than $L_0$ for all sufficiently large $N$. 
In other words,  there will be an infinitely large region as $N$ tends to infinity, namely $1\leq L \leq  T_{\BC:\HQ}$, where $L_0$ falls into and where BC penalizes more than HQ.
As a result, asymptotically the bridge criterion does not overfit.
On the other hand, the bridge criterion will not underfit because the largest penalty preventing from selecting $L+1$ versus $L$ is $2\LM /N$, which will be  less than any fixed positive number $g_{L_0}$ defined in (\ref{gain_2}) for all sufficiently large $N$. The bridge criterion is therefore consistent.

The inequality $(2\LM / N) /k \leq 2/N$ for any $1 \leq k \leq \LM$  guarantees that BC penalizes more than AIC so that it does not cause much overfitting even in the case of small $N$ or large $L_0$.
Since BC penalizes less for larger orders and finally becomes similar to AIC, it is able to share the asymptotic optimality of AIC under suitable conditions.
To further illustrate why the bridge criterion is expected to work well in general, we make the following intuitive argument about the model selection procedure.
As we shall see, the bent curve of BC well connects BIC (or HQ) and AIC so that a good balance between the underfitting and overfitting risks is achieved. 
The rigorous theory will be established in Section~\ref{sec:theory}.


\textbf{Intuitive argument:}

\begin{figure}[h!]
\centering
  \includegraphics[width=0.5\linewidth]{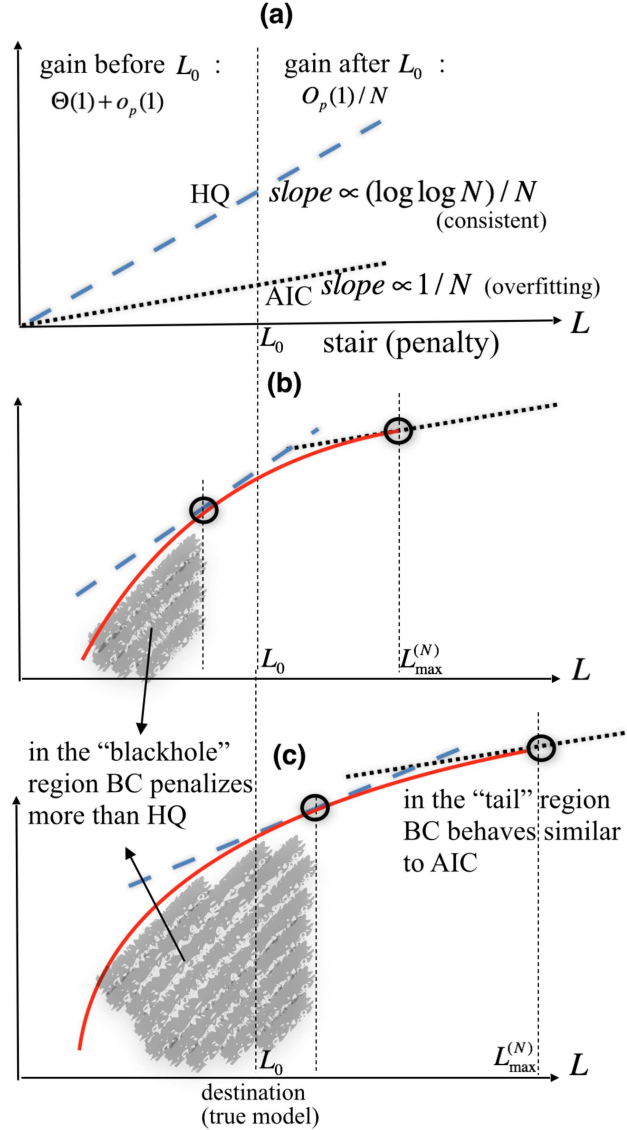}
  \vspace{-0.0in}
  \caption{(a) Curve $J_{\AIC}$ (blue dash) and $J_{\HQ}$ (black small dash), (b) the joint plot of $J_{\BC}$ (red thick line) and $J_{\AIC}, J_{\HQ}$, by shifting the latter two to be tangent to $J_{\BC}$ at tangent points $T_{\BC:\AIC}, \ T_{\BC:\HQ}$ (circled), in which $T_{\BC:\AIC} < L_0$, and (c) the evolution of plot (b) to the scenario  $T_{\BC:\HQ} \geq L_0$ as $N$ increases}
  \label{fig:tale}
  \vspace{-0.0in}
\end{figure}

To gain further intuition, we consider
an insect who is climbing a slope that is determined by a particular penalty curve $J(L)$ from the starting point $L=1$ to the maximal possible end $m=\LM$ (Fig.~\ref{fig:tale}).  
Fig.~\ref{fig:tale}(a) illustrates $J_{\AIC}(L)$ (black small dash) and $J_{\HQ}(L)$ (blue dash).
We only drew $J_{\HQ}(L)$ for brevity, as there is no essential difference between the two strongly consistent criteria HQ and BIC.  

\textit{The climbing scheme and the goal: }
At each step $L$, the insect moves to step $L+1$ if its gain is larger than its loss, and it will not move any more once it stops. The gain refers to the increased goodness of fit to the data (which is $\hat{g}_{L+1}$ in our autoregressive model), the loss refers to the  penalty of increased model complexity (which is $J(L+1)-J(L)$), and the last step where the insect stops is denoted by $\hat{L}$.
The goal is to design a proper slope such that the insect stops at a ``desired destination'' that will be elaborated on below.

\textit{The tangent points of two slopes: } 
A slope can be written as $\sum_{k=1}^L t_k$. The performance of the insect is determined by each increment $t_k$, and is not affected if the slope is shifted by any constant that does not depend on $L$.
We thus shift the curves $J_{\AIC}(L)$ and $J_{\HQ}(L)$ to be tangent to the log-like curve of $J_{\BC}(L)$.  
By our design of $J_{\BC}(L)$, the tangent points between $J_{\BC}(L)$ and $J_{\AIC}(L),J_{\HQ}(L)$ curves are respectively at steps $T_{\BC:\AIC} = \LM$, $T_{\BC:\HQ} = 2\LM/ (c\log \log N)$.
Before step $T_{\BC:\HQ}$, the insect on BC slope suffers more loss than on HQ slope in each move, while the other way around after step $T_{\BC:\HQ}$.

\textit{The well-specified case: }
Now we categorize two distinct scenarios: where the desired destination is within finitely many steps, and where the desired destination is beyond finitely many steps.
In the former case, there is a clear target step $L_0$. A good slope should be designed such that the insect stops at step $L_0$.  
It is already known in this case that HQ slope is good while AIC slope is not. In fact, it can be illustrated by Fig.~\ref{fig:tale}(a), in which the gain after $L_0$ is $O_p(1)/N$, 
smaller than $\Theta(\log\log N) / N $ (which is usually guaranteed by the law of the iterated logarithm) while larger than $O(1)/N$ with a positive probability for sufficiently large $N$. 
How about BC? 
It is worth mentioning that our argument for the insect is implicitly built upon $N$, and the concept of consistency is about large $N$ asymptotics. Suppose that $N$ keeps increasing, the aforementioned tangent step $T_{\BC:\HQ}$  will be not only larger than $L_0$ but also diverging to infinity given that $\log \log N = o(\LM)$.
In other words, there is the ``blackhole'' region $[0,T_{\BC:\HQ}]$ (Fig.~\ref{fig:tale}(b) and (c)), in which BC slope is steeper than HQ slope, and which grows to be infinitely large. It results in two consequences: 
First, the insect will find it more and more difficult to escape from the region because the increased loss from moving each step needs to be compensated by its gain. Take the autoregressive models as an example. After moving each step the gain is approximated independent $\chi_1^2/N$, the expectation of which is less than the loss $2/N$; so the probability of the cumulated sum of gains being larger than that of loss decreases to zero rapidly as the number of steps increases. 
Second, once the insect is trapped in the blackhole, it encounters more difficulty to move forward on a BC slope than on a HQ one. Since on the HQ slope the insect will not move beyond step $L_0$ (due to the strong consistency of HQ), on a BC slope it will not, either. 

On the other hand, the insect will not stop before step $L_0$. This is because of two facts: First, the largest penalty preventing from moving forward is $J_{\BC}(1)=o(1)$; Second, the gain of the insect moving from step $L$ to $L+1$ when $L < M_0$ is usually at least $\Theta(1)+o_p(1)$ (which is true when $\psi_{L+1,L+1}\neq 0$ in autoregressive models). 
Therefore, the insect stops at step $L_0$ on a BC slope.

\textit{The mis-specified case: }
The fact that $T_{\BC:\AIC}=\LM$  guarantees that BC slope is always steeper than AIC slope so that the insect does not move too far.
Because the BC slope is in a concave shape, the insect moves easier and easier for larger steps. In the case where the appropriate destination tends to infinity, the insect will soon move to the tail part of the slope. 
As one can see from Fig.~\ref{fig:tale}(c), in the tail part the slope is designed to be similar to AIC (and it becomes exactly AIC at the end step $L=\LM$), it is possible to share the asymptotic optimality of AIC.

In summary, the bent curve of the BC well connects AIC and HQ so that a good balance between the underfitting and overfitting risks can be achieved.
We emphasize that the above argument does not match exactly to the rigorous proof, since the decision making of the insect is carried out sequentially, while the aforementioned criteria select $\hat{L}$ via global optimum.
Nevertheless, the argument for the insect does shed some light on why BC is likely to perform in the way we desire: to automatically behave like a consistent one while the underlying model is well-specified, and an efficient one otherwise, alleviating the risk caused by an analyst's initial prejudice. 
Besides this, the above argument does not assume any concrete probabilistic model, and thus it seems to be a promising criterion for other statistical inference as well.

\section{Performance of the bridge criterion} \label{sec:theory}

In this section, we establish rigorous theory on the consistency and efficiency of the bridge criterion proposed in (\ref{TS_criterion}). 
We prove its consistency and asymptotic efficiency in Subsection~\ref{sec:Consistency} and \ref{sec:efficiency}, respectively.
In Subsection~\ref{sec:modifiedBC}, we propose an extended bridge criterion and its associated two-step strategy, in order to relax some technical assumptions.
In view of the above intuitive argument, the extended  criterion works in the following way. 
Let the insect clime on the AIC slope,  and record its ending point $\Laic$; modify the BC increment $J_{\BC}(L)-J_{\BC}(L-1)$ from $2\LM /(NL)$ to $2\M /(NL)$, where $\M$ is slightly smaller than $\LM$; let the insect move again on the modified BC slope with boundary $\Laic$. 
In this way, the insect can still stop at $L_0$ if it is finite, and otherwise moves faster towards the end $\Laic$ as if it were on the AIC slope. 

\subsection{Consistency}
\label{sec:Consistency}

%
\begin{theorem} \label{thm:consistency}
Suppose that the time series data is generated from a finite order autoregression, and that 
\begin{align}\label{condition_1}
\lim\limits_{N \rightarrow \infty} \frac{\LM }{ \log \log N } = \infty , \quad \lim\limits_{N \rightarrow \infty}  \frac{\LM }{ N^{\frac{1}{2}} }   = 0 \ .
\end{align}
Then the bridge criterion is consistent.
In addition, if $\hat{L}$ is selected from any finite set of integers that does not depend on $N$ and that contains the true order $L_0$, then $\hat{L}$ converges not only in probability but also almost surely to $L_0$. 
\end{theorem}

\begin{remark} \label{remark_new}
\normalfont
Theorem~\ref{thm:consistency} proves the consistency of bridge criterion under mild assumptions.
It is worth mentioning that 
an analyst does not suffer the risk of specifying a finite candidate set $\{1,\ldots,\lm\}$ that excludes the true order $L_0$.
Because any finite true order will be eventually included as a candidate and evaluated by bridge criterion, as the sample size becomes large.
The  proof of Theorem~\ref{thm:consistency} is given in the Supplementary Material.

\end{remark}

\begin{remark}
\normalfont
The proof of Theorem~\ref{thm:consistency} could be adapted in such a way that $\lim_{N \rightarrow \infty} \LM / \log \log N$ $= \infty$ is not necessary to prove the consistency. It was for proving the strong consistency of  $\hat{L}$  if  any finite candidate set including the true order $L_0$ is specified instead. 
Nevertheless, various numeric experiments show that this condition greatly enhances the performance of the bridge criterion under finite sample size when applied to the candidate set $\{1,\ldots,\LM\}$. 
\end{remark}

\subsection{Asymptotic efficiency}
\label{sec:efficiency}

We  introduce the following  notation.
The matrix norm $\norm{\cdot}$ is defined by $\norm{M} = \sup_{\norm{y}_2=1} \norm{ M y }_2$, where $\norm{\cdot}_2 $ denotes the Euclidean norm of a column vector.
For a positive definite matrix $A$, the norm $\norm{\cdot}_{A}$ is defined by 
$\norm{y}_{A} = (y^\T A y )^{1/2}.$
If two vectors $y_1=[y_{1,1},\ldots,y_{1,L_1}]^\T$ and $y_2=[y_{2,1},\ldots,y_{2,L_2}]^\T$ are of different sizes, then
we allow subtraction of those vectors by modifying the definition in the following way.
Given $y_1,y_2$, define $y_1' , y_2' $ as vectors of size $L' = \max\{L_1, L_2\}$ by appending $\max\{L_1, L_2\}-\min\{L_1, L_2\}$ zeros to the tail of $y_1$ or $ y_2$.
We define subtraction of $y_1,y_2$ in this case as $y_1' - y_2'$.
Similarly, if
the size of a vector $y$ is smaller than a positive definite matrix $A$ of size $k \times k$, $\norm{y}_A$ is the same as $\norm{y'}_A$ where $y'$ is of size $k$ by appending zeros to the tail of $y$.

We are usually interested in the one-step prediction error if a mismatch filter, as defined below, is specified  \cite{akaike1969fitting,akaike1970statistical}.
Assume that the data is generated from a filter $ \Psi_{L_0}$ as in  (\ref{AR}).
The one-step prediction error of using filter $ \Lambda_{L}$  minus that  of using the true filter is referred to as mismatch error
\begin{align} \label{mismatch}
E  \bigl\{ [x_n,\ldots,x_{n-L'+1}] ( \Psi_{L_0} -  \Lambda_{L}) \bigr\}^2
= \norm{\Lambda_{L} - \Psi_{L_0}}_{\Gamma_{ L'}}^2 \ . 
\end{align}
where $ L' = \max\{L_0,L\}$ and $\Gamma_{ L'}$ is the $ L' \times  L'$ covariance matrix of the true autoregression, namely its $(i,j)$th element is $\gamma_{i-j}$. 
%
%
The following assumptions are  needed  for this section.
\begin{assumption} \label{A1}
The data $\{x_n: n=1,\ldots, N_0\}$  
is generated from the recursion
    $
         x_n + \psi_{\infty,1} x_{n-1} +  \psi_{\infty,2} x_{n-2}  + \cdots   =\varepsilon_n,
    $ 
where $\psi_{\infty,j} \in \mathbb{R}$, $\sum_{j=1}^{\infty} | \psi_{\infty,j} | < \infty $,  $\varepsilon_n$'s are independent and identically distributed according to $\mathcal{N}(0,\sigma^2)$, and
the associated power series $\Psi(z) = 1 + \psi_{\infty,1} z^{-1} +\psi_{\infty,2} z^{-2} + \cdots$ converges and is not zero for $|z| \geq 1$.
\end{assumption}


\begin{assumption} \label{A2}
$\{\LM \}$ is a sequence of positive integers such that $\LM \rightarrow \infty$ and $\LM =o( N^{1/2} ) $ as $N$ tends to infinity.
\end{assumption}

\begin{assumption} \label{A3}
The order of the autoregressive process (or the size  of filter $ \Psi_{\infty}$) is infinite.
\end{assumption}

\begin{remark} \label{remark:ass}
\normalfont
Assumption~\ref{A1} is a more general assumption than we had in previous sections. 
Under Assumption~\ref{A1}, we have
$0 < \gamma_0 = \norm{\Gamma_1 } \leq \norm{\Gamma_2 } \leq \cdots \leq \norm{\Gamma }$, where $\Gamma =[ \gamma_{i-j}]_{i,j=1}^{\infty}$
 is the infinite dimensional covariance matrix with norm
 $
 \norm{\Gamma} = \sup_{\norm{y}_2 = 1} \left[ \sum_{i=1}^{\infty} \{\sum_{j=1}^{\infty} \gamma_{i-j} y_j \}^2 \right]^{1/2}.
 $

Assumption~\ref{A3} has been assumed in several technical lemmas in \cite{shibata1980asymptotically} that we are going to introduce. 
For those lemmas and the scope of this paper, Assumption~\ref{A3} can be generalized to allow for the case where the order of the autoregressive process, denoted by $L_0(N)$, is finite but depends on $N$. In other words, the data generating process varies with $N$. In that case, the associated power series that appeared in Assumption~\ref{A1} may be written as $\Psi_{N}(z)=1 + \psi_{L_0(N),1} z^{-1} +\cdots+ \psi_{L_0(N),L_0(N)} z^{-L_0(N)} $, and that assumption is accordingly replaced by: $\Psi_{N}(z)$ is not zero for $|z| \geq 1$ and it converges as $N$ tends to infinity; an additional requirement is the divergence of $L_N^{*}$ (introduced below) as $N$ tends to infinity.
\end{remark}

In this section, we show that the proposed order selection criterion asymptotically minimizes the mismatch error under certain conditions.
Define the cost function $C_N(L) = L \sigma^2/N +  \norm{ \Psi_L - \Psi_{\infty}}_{\Gamma}^2   $. 
It can be regarded as the expected mismatch error if an estimated filter of order  $L$ is used for prediction.
%
%
In fact, under Assumptions~\ref{A1}--\ref{A3}, it holds that \cite[Proposition 3.2]{shibata1980asymptotically}
\begin{align} \label{lemma_prop42}
\lim\limits_{N \rightarrow \infty} \max_{1 \leq L \leq \LM} \abs[\bigg]{ \frac{ \norm{\hat{ \Psi}_{L}  - \Psi_{\infty}}_{\Gamma}^2  }{C_N(L) } - 1 } = 0 \quad \textrm{ in probability}.
\end{align}
%

In addition, if we use $\{L_N^{*}\}$ to denote a sequence of positive integers which achieves the minimum of $C_N(L)$ for each $N$, namely
$L_N^{*} = \argmin_{1 \leq L \leq \LM} C_N(L) $,
then for any random variable $\tilde{L}$ possibly depending on $ \{x_n: n=1,\ldots,N\}$, and for any $\epsilon >0$,
it holds that  
$\lim_{N \rightarrow \infty} \P  \bigl\{  \norm{ \hat{ \Psi}_{\tilde{L}} - \Psi_{\infty}}_{\Gamma}^2  / C_N(L_N^{*})  \geq 1-\epsilon \bigr\} = 1$ \cite[Theorem 3.2]{shibata1980asymptotically}.
The result shows that the cost of the estimate $ \hat{ \Psi}(\tilde{L})$ is no less than $C_N(L_N^{*})$ in probability for any order selection $\tilde{L}$.
An order selection $\tilde{L}$ is called asymptotically efficient if
\begin{align} \label{prop42_equivalence}
\lim\limits_{N \rightarrow \infty} \frac{  \norm{ \hat{ \Psi}_{\tilde{L}} - \Psi_{\infty}}_{\Gamma}^2  }{C_N(L_N^{*}) } = 1  \quad \textrm{ in probability}.
\end{align}
Equality (\ref{prop42_equivalence}) can be equivalently written as 
$\lim_{N \rightarrow \infty}  C_N(\tilde{L}) /C_N(L_N^{*})  = 1  \textrm{ in probability} $
in view of Equality~(\ref{lemma_prop42}).
%
The following result establishes the asymptotic efficiency of bridge criterion in two common scenarios, i.e., where the mismatch error $ \norm{ \Psi_{L} - \Psi_{\infty}}_{\Gamma}^2 $ decays algebraically or exponentially in $L$.
The two cases cover a wide range of linear processes as we point out in Remark~\ref{remark1}.  
Its proof is given in the Supplementary Material. 

\begin{proposition} \label{thm:efficiency}
Suppose that Assumptions~\ref{A1}--\ref{A3} hold.
\begin{enumerate}
\item
Suppose that the mismatch error $ \norm{ \Psi_{L} - \Psi_{\infty}}_{\Gamma}^2 $ satisfies
\begin{align} \label{con1_1}
\log  \norm{ \Psi_{L} - \Psi_{\infty}}_{\Gamma}^2 = -\gamma \log L + \log c_L
\end{align}
where $\gamma \geq 1$ is a constant, and the series $\{c_L: L=1,2, \ldots\}$ is lower bounded by a positive constant and $c_{L+1} / c_{L} < 1+ \gamma /(L+1) $.
If
\begin{align} \label{case1}
\LM = O\biggl( N^{\frac{1}{1+\gamma}-\varepsilon} \biggr)
\end{align}
holds for a fixed constant $0 < \v < 1/(1+\gamma)$, then the bridge order selection criterion is asymptotically efficient.

\item
Suppose that the mismatch error satisfies the equality  
\begin{align} \label{con2_2}
\log  \norm{ \Psi_{L} - \Psi_{\infty}}_{\Gamma}^2 = -\gamma  L + \log c_L
\end{align}
where $\gamma > 0$ is a constant, and the series $\{c_L: L=1,2, \ldots\}$ is lower bounded by a positive constant and $c_{L+1} / c_{L} \leq  q$ for some constant $q< \exp(\gamma)$.
If
\begin{align} \label{case2}
\LM \leq \frac{1-\varepsilon}{\gamma} \log N  
\end{align}
holds for a fixed constant $0 < \varepsilon < 1$, then the bridge order selection criterion is asymptotically efficient.

\end{enumerate}

\end{proposition}

%
%
%
%

\begin{remark} \label{remark1}
\normalfont
To provide an intuition of condition (\ref{con1_1}), in view of Remark~\ref{remark:ass} we prove that if the order of autoregressive process is not infinity but $L_0(N)$ (which grows with $N$) instead, and if $ \Psi_{L_0(N)}$ is uniformly distributed in 
$S_{L_0(N)}(1)$ for any given $N$, 
then for large $L$ ($1 \leq L \leq L_0(N)$)
\begin{align} \label{gamma_1}
E \bigl\{ \log  \norm{ \Psi_{L} - \Psi_{L_0(N)}}_{\Gamma_{L_0(N)}}^2 \bigr\} = - \log L + \log L_0(N) +o_L(1) \ .
\end{align}
The proof is given in the Supplementary Material.
Furthermore, it is known that  condition (\ref{con2_2}) holds (with constant series $c_L$) when the data is generated from a finite order moving-average process \cite{shibata1980asymptotically}. 
\end{remark}


However, the proposed bridge criterion in (\ref{TS_criterion}) is not fully satisfactory in terms of asymptotic efficiency.
For BC to achieve efficiency, our Proposition~\ref{thm:efficiency} requires $ \LM$ to satisfy 
(\ref{case1}) or (\ref{case2}) depending on the  underlying mismatch error.  
This poses two concerns: first, the mismatch error as a function of $L$ is usually unknown in advance, and it can be more complex than those characterized by (\ref{con1_1}) and (\ref{con2_2}); 
second, the chosen $\LM$ is not large enough to incorporate all possible competitive models into the candidate set; this is because $ \LM$ is always $\v$-away (in terms of the order) to the minimum of $C_N(L) $ over all positive integers $L \in \mathbb{N}$.  
This has motivated us to extend the bridge criterion in such a way that 1) it relaxes the conditions required by (\ref{con1_1}) and (\ref{con2_2}), and 2) it selects the optimal order from a broad candidate set, 
and 3) it still achieves either consistency in the well-specified case or efficiency in mis-specified case.

\subsection{Adaptive selection of $\LM$}
\label{sec:modifiedBC}

To achieve the aforementioned goal, we propose a general strategy that consists of two steps. 
\begin{enumerate}[1.]
\item choose any $\LM = o(\sqrt{N})$ and apply AIC to obtain 
$\Laic$; 
\item within the range $1,2,\ldots,\Laic$, select the optimal order (denoted by $\Lbc$) by minimizing the modified BC penalty  
\begin{align} \label{BC2}
\BC(N,L) = \frac{2\M }{N} \sum\limits_{k=1}^L \frac{1}{k} 
\end{align}
where $\M$ is a number to be chosen. 
\end{enumerate}

We note that $\M = \LM$ was chosen in the previous sections, but it may not be the ideal choice in our two-stage approach, as we shall see later. 
We define 
\begin{align} \label{universal_opt}
  \Lopt = \argmin_{L \in \mathbb{N}} C_N(L).
\end{align}
to be the ``universally optimal order''. 
In most cases $\Lopt$ is upper bounded by $N^{1-\v}$ for a fixed $\v > 0$. 
For instance, if $\norm{\psi_L - \psi_{\infty}}_{\Gamma}^2$ follows the algebraic decay 
$cL^{-\gamma}$ for some $\gamma > 0$, then $\Lopt = \Theta \bigl(N^{1/(1+\gamma)} \bigr)$.
Nevertheless, it is possible that $\Lopt > \sqrt{N}$.

In the rest of this section, we consider the case $\Lopt \leq \LM$ in order to: \\
1) take into account the most competitive model that does not depend on the choice of $\LM$, as (\ref{universal_opt}) implies
$
  \Lopt = \argmin_{1 \leq L \leq \LM} C_N(L)
$; \\ 
2) simplify technical derivations. But we emphasize that this requirement is not essential.
\begin{assumption} \label{A4}
In the mis-specified scenario, it holds that $\Lopt \leq \LM$.
In addition, $C_N(L)$ has a well-separated mode in the sense that 
if $\lim_{N\rightarrow \infty} C_N(L_N)/C_N(\Lopt)=1$ holds for a sequence $L_N$, then $\lim_{N\rightarrow \infty} L_N/L(\Lopt)=1$.
\end{assumption}
\begin{remark} \label{remark:ass4}
\normalfont
The efficiency of AIC under mis-specified model implies that \\
$\lim_{N\rightarrow \infty} C_N(\Laic)/C_N(\Lopt) = 1$
in probability
which, given Assumption~\ref{A4}, further implies 
\begin{align}
	&\lim\limits_{N\rightarrow \infty} \frac{\Laic}{\Lopt} = 1 
	\quad \textrm{in probability.} \label{new_condition2}
\end{align} 
Assumption~\ref{A4} is easily satisfied in many common cases. 
For example, we consider two common scenarios that were also described in Proposition~\ref{thm:efficiency}:
the mismatch error has an algebraic decay 
$\norm{ \Psi_{L} - \Psi_{\infty}}_{\Gamma}^2 = c L^{-\gamma}$, or an exponential decay 
$\norm{ \Psi_{L} - \Psi_{\infty}}_{\Gamma}^2 = c\exp(-\gamma  L) $.
We let $q_N = L_N/\Lopt$. 
Via straightforward calculation, $\lim_{N\rightarrow \infty} C_N(L_N)/C_N(\Lopt)=1$ can be rewritten as 
$
\lim_{N\rightarrow \infty}  q_N^{-\gamma}(1 + \gamma q_N^{\gamma+1} / (1+\gamma ) = 1 
$ 
in the case of algebraic decay, 
and it can be rewritten as
$
\lim_{N\rightarrow \infty}  \exp\{-\gamma(L_N-\Lopt)\}/(1+\gamma \Lopt) + \gamma L_N / (1+ \gamma \Lopt) = 1
$ 
in the case of exponential decay.
In both cases, it follows that 
$\lim_{N\rightarrow \infty} q_N = 1$ in probability.
\end{remark}

The following theorem establishes the consistency and efficiency of the  two-stage strategy.

\begin{theorem} \label{thm:modifiedBC}
Suppose that $\Laic$ is obtained from the first step of the two-step strategy, and Assumption~\ref{A2} holds.
Suppose that there exists   a sequence $\M$ (indexed by $N$) satisfying 
\begin{align}
	&\lim\limits_{N \rightarrow \infty} \frac{\M}{\log \log N} = \infty.\label{new_condition5} 
\end{align}
In addition, assume that under a mis-specified model class, Assumptions~\ref{A1},\ref{A3},\ref{A4} hold, and for all sufficiently large $N$
\begin{align}
&\M \leq \frac{q\Lopt}{\log \Lopt}  ,
\label{new_condition3} 
\end{align}
where  $0 < q < 1$ is some constant and $\Lopt$ was defined in (\ref{universal_opt}).
Then, using the above two-stage strategy, the modified bridge criterion in (\ref{BC2}) is consistent in the well-specified case and efficient in the mis-specified case. Moreover, if in the well-specified case $\hat{L}$ is selected from a finite candidate set that does not depend on $N$ and that contains the true order $L_0$, then $\hat{L}$ converges almost surely to $L_0$. 
\end{theorem}

%
\begin{remark} \label{M_choiceInPractice}
\normalfont
  We note that Conditions (\ref{new_condition5}) and (\ref{new_condition3}) are fairly weak. For instance, $\Lopt$ is respectively $\Theta(N^{r})$ ($0<r<1$) and $\Theta(\log N)$ in the two cases described in Proposition~\ref{thm:efficiency}, so we may choose $M_N=(\log N)^{\tau}$ with any $0<\tau<1$. 
\end{remark}

\begin{remark}
\normalfont
We provide an intuitive reasoning here. In the well-specified scenario, 
	(\ref{new_condition5}) guarantees consistency due to Theorem~\ref{thm:consistency}. 
In the mis-specified scenario, (\ref{new_condition2}) and (\ref{new_condition3}) imply that $\M < \Laic / \log \Laic $. Such $\M$ produces penalty increments $J_{\BC}(L+1)-J_{\BC}(L)$ that are lighter than AIC for large $L$ (recall that the candidate set in the second step is $1,\ldots,\Laic$). In view of that, BC produces $\hat{L}$ that is close to the boundary $\Laic$. 
\end{remark}

\begin{remark}
\normalfont
Another form of the modified bridge criterion is written as 
\begin{align} \label{BC3}
\BC(N,L) = \frac{2\M }{N} \sum\limits_{k=1}^L k^{-\zeta}
\end{align}
where $\zeta>0,\zeta \neq 1$.
By a similar proof, it can be shown that Theorem~\ref{thm:modifiedBC} can be modified to the case $0< \zeta < 1$ by requiring the following changes: 
	replace $ \Lopt / \log \Lopt $ by $(\Lopt)^{\zeta}$ in (\ref{new_condition3}), and require $q < 1-\zeta$. 
In addition, Theorem~\ref{thm:modifiedBC} can be modified to the case $\zeta > 1$ via replacing 
	$ \Lopt / \log \Lopt $ by $\Lopt/a(\zeta)$ in (\ref{new_condition3}),
where $a(\zeta) = \sum_{k=1}^{\infty} k^{-\zeta}$. 
As a possible future work, it would be interesting to compare the performance of $\zeta=1$ and $\zeta \neq 1$.


%
%
%
%
%
%
%
%
%
%

\end{remark}

\begin{remark} \label{remark:PI} 
\normalfont
	Building upon the proposed bridge criterion, we define the following parametricness index (PI):
	\begin{align}\label{PI}
	\PI = \left\{
	\begin{aligned} 
		&\frac{|\Lbc - \Laic|}{|\Lbc - \Laic| + |\Lbc - \Lbic|} & \textrm{ if }\Laic \neq \Lbic  \\
		& 1 & \textrm{ otherwise. } 
	\end{aligned}
	\right. 
	\end{align}
	Following the definition, $\PI \in [0,1]$. Intuitively, $\PI$ is close to one in the well-specified model class where $\Lbc,\Lbic$ do not differ much, while close to zero in a mis-specified one where $\Lbc,\Laic$ are close and much larger than $\Lbic$.  
	The goal of PI is to measure the extent to which the specified model class is adequate in explaining the observed data, namely to assess the confidence that the selected model can be practically treated as the data-generating model. The larger $\PI$, the more confidence. 
	Similar concept has been introduced in \cite{liu2011parametric} for the goal of estimating the regression function.	
	The following proposition shows that 
	$\PI$ converges in probability to one for the  well-specified case. Though we cannot prove that $\PI$ converges in probability to zero for various mis-specified cases in general, for illustration purpose we prove for some typical mis-specified cases. 
	Experiments on various synthetic data in Section~\ref{sec:Num_Results} have shown that  $\PI$ performs in the way we expected. 
\end{remark}

\begin{proposition} \label{prop:PI}
	Under the same conditions of Theorem~\ref{thm:modifiedBC}, if the  model class is well-specified, $\PI$ converges in probability to one as $N$ goes to infinity;	
	If  the  model class is mis-specified, and we further assume that $C_N(L) + (\log N - 2)L \sigma^2/N $ achieves its minimum at $L_{*}^{(N)}$ and  $ \lim_{N \rightarrow \infty} L_{*}^{(N)} / \Lopt = 0$, then $\PI$ converges in probability to zero as $N$ goes to infinity.
	An example is where the  
	 mismatch error satisfies 
	 $\norm{ \Psi_{L} - \Psi_{\infty}}_{\Gamma}^2 = c L^{-\gamma }$, 
where $\gamma$ and $c$ are positive constants.
\end{proposition}

\section{Numerical results}
\label{sec:Num_Results}

In this section, we present experimental results to demonstrate the theoretical results and the advantages of bridge criterion on both synthetic and real-world datasets. 
Throughout the experiments, we use the two-step bridge criterion defined in (\ref{BC2}), and we adopt 
\begin{align} \label{eq:newer11}
\LM = \lfloor N^{1/3} \rfloor, \quad \M = (\log N)^{0.9}
\end{align} 
due to Theorem~\ref{thm:modifiedBC} and Remark~\ref{M_choiceInPractice}, where $N$ is the sample size.   

\subsection{Synthetic data experiment: consistency in finitely dimensional model}

The purpose of this experiment is to show the consistency of BC and BIC. 
The performance of BC, AIC, and BIC in terms of order selection for well-specified model class is summarized in Table~\ref{table:result1}.
In Table~\ref{table:result1}, the data is simulated using autoregressive filters $ \Psi_2 = [\alpha, \alpha^2]^{\T}$ for $\alpha=0.3,-0.3,0.8,-0.8$. 
For each $\alpha$, the estimated  orders are tabulated for $1000$ independent realizations of AR$(2)$ processes
$x_n+\alpha x_{n-1}+\alpha^2 x_{n-2}=\epsilon_n, \, \epsilon_n \sim \mathcal{N}(0,1)$.
The experiment is repeated for different sample sizes $N=100,500,1000,10000$. 
As was expected, the performance of the bridge criterion lies in between AIC and BIC, and it is consistent  when $N$ tends to infinity. In addition, the convergence for $\alpha=0.3,-0.3$ is slightly slower compared with $\alpha=0.8,-0.8$,  because of their smaller signal to noise ratios. 

\begin{table}[h!]
    \caption{\label{table:result1} Selected orders for AR$(2)$ processes (1000 realizations for each $\alpha$ and $N$)}
\centering
\fbox{
    \begin{tabular}{ c  c | c c c | c c c | c c c | c c c }
    \multicolumn{2}{c|}{ }	&\multicolumn{3}{c|}{$N=$100} &\multicolumn{3}{c|}{$N=$500} &\multicolumn{3}{c|}{$N=$1000} &\multicolumn{3}{c}{$N$=10000} \\
    \multirow{1}{1.5em}{$\alpha$} &$\hat{L}$ &BC &AIC  &BIC  &BC &AIC  &BIC   &BC &AIC &BIC   &BC &AIC &BIC 
    \\
	\hline
	\multirow{4}{1.5em}{$0.3$} 	&1	 &784 &548 &851   &558 &213 &661	&298	 &51	 &405   &0	     &0	 &0  
	\\
	&\cellcolor{lightgray}2 	 &\cellcolor{lightgray}151 &\cellcolor{lightgray}292 &\cellcolor{lightgray}135	
	&\cellcolor{lightgray} 372 &\cellcolor{lightgray}558 &\cellcolor{lightgray}333
	&\cellcolor{lightgray}619	 &\cellcolor{lightgray}677 &\cellcolor{lightgray}589  
	&\cellcolor{lightgray}949  &\cellcolor{lightgray}720	 &\cellcolor{lightgray}999  
	\\
	&3	 &36	  &98  &13	 &37 &113 &5  	&38	 &125  &5	    &21    &97	 &1 
	\\
	&$>3$	 &29	 &62	  &1	   &33 &116 &1     &45	 &147 &1	    &30	     &183	 &0  
	\\
	\hline
	\multirow{4}{1.5em}{$-0.3$}     	&1	 &777 &566 &845  &535 &208 &628   &297	 &45	 &375   &0	     &0	 &0  
	\\
	&\cellcolor{lightgray}2	 &\cellcolor{lightgray}166 &\cellcolor{lightgray}301 &\cellcolor{lightgray}145 
	&\cellcolor{lightgray}392 &\cellcolor{lightgray}536 &\cellcolor{lightgray}365	
	&\cellcolor{lightgray}624	 &\cellcolor{lightgray}688 &\cellcolor{lightgray}617  
	&\cellcolor{lightgray}958  &\cellcolor{lightgray}719	 &\cellcolor{lightgray}997  
	\\
	&3	 &28	  &64  &8	   &32 &110 &6   	&32	 &112  &7	    &22    &122	 &3 
	\\
	&$>3$	 &29	 &69	  &2	 &41 &146 &1     	&47	 &155 &1	    &20	     &159	 &0  
	\\
	\hline
	\multirow{4}{1.5em}{$0.8$}   &1	 &0   &0 &0       &0  &0  &0       &0  &0  &0             &0   &0	 &0  
	\\
	&\cellcolor{lightgray}2	 &\cellcolor{lightgray}823 &\cellcolor{lightgray}749 &\cellcolor{lightgray}957 
	&\cellcolor{lightgray}891 &\cellcolor{lightgray}734 &\cellcolor{lightgray}988  
	&\cellcolor{lightgray}906	 &\cellcolor{lightgray}715 &\cellcolor{lightgray}992  
	&\cellcolor{lightgray}944  &\cellcolor{lightgray}726	 &\cellcolor{lightgray}998  
	\\
	&3	 &102	  &148  &36   &44 &125 &11  &41	 &118  &8	    &24    &102	 &2 
	\\
	&$>3$	 &75	 &103	  &7	&65 &141 &1      &53	 &167 &0	    &32  &172	 &0  
	\\
	\hline
	\multirow{4}{1.5em}{$-0.8$} &1	 &0   &0 &0    &0  &0  &0      &0  &0  &0             &0   &0	 &0  
	\\
	&\cellcolor{lightgray}2	 &\cellcolor{lightgray}860 &\cellcolor{lightgray}783 &\cellcolor{lightgray}968
	&\cellcolor{lightgray}876 &\cellcolor{lightgray}738 &\cellcolor{lightgray}980 
	&\cellcolor{lightgray}878	 &\cellcolor{lightgray}709 &\cellcolor{lightgray}994  
	&\cellcolor{lightgray}949  &\cellcolor{lightgray}703	 &\cellcolor{lightgray}999  
	\\
	&3	 &82	  &127  &29  &54 &112 &18  &55	 &133  &5	    &23    &115	 &1 
	\\
	&$>3$	 &58	 &90   &3	 &70 &150 &2    &67	 &158 &1	    &28  &182	 &0  
	\\
    \end{tabular}}
\end{table}

\subsection{Synthetic data experiment: efficiency in finitely and infinitely dimensional models}

The purpose of this experiment is to show that the proposed order selection criterion achieves the asymptotic efficiency for both the well-specified and the mis-specified cases.
The performance of BC in terms of mismatch error is compared with those of AIC and BIC in Table~\ref{table:result2}.
Recall that the mismatch error defined in (\ref{mismatch}) is the expected one-step ahead prediction error minus the variance of noise, when an estimated filter is applied to an independent and identically generated dataset. 
We consider three different data generating processes below.
In Table~\ref{table:result2}, for each case and sample size $N=100,500,1000,10000$, the tabulated mismatch error produced by each criteria were the mean of 1000 repeated independent experiments. 
The mean parametricness index defined in Remark~\ref{remark:PI} (denoted by $\PI$) in each case was also tabulated. 

\textbf{Case 1: }
The first case is AR$(1)$ with $ \Psi_1 = [0.9]$, namely $x_n+0.9 x_{n-1}=\epsilon_n, \, \epsilon_n \sim \mathcal{N}(0,1)$. This is a well-specified model. As we can see, once the true order is selected with probability close to one, the resulting predictive performance is also asymptotically optimal.

Here, we briefly explain how to calculate the exact mismatch error in (\ref{mismatch}) for any estimated filter  of size $L$ that may or may not equal to $L_0$. If suffices to express the covariance matrix $\Gamma_{ L'}$ or its elements $\gamma_0,\ldots,\gamma_{L'-1}$ in terms of the known $\Psi_{L_0}$, where $L' = \max \{L_0, L\} $. 
We define the correlation vector and matrix by
$\rho_{L_0} = [\gamma_1/ \gamma_0, \ldots, \gamma_{L_0}/\gamma_0]^\T , P_{L_0} =   \Gamma_{L_0}/ \gamma_0 $.
By rewriting the Yule-Walker equation $P_{L_0} \Psi_{L_0} = -\rho_{L_0} $, we obtain 
$(I + \Phi_{L_0} ) \ \rho_{L_0} = - \Psi_{L_0} $ where 
{\small
\begin{align*}
 \Phi_{L_0} =
\begin{bmatrix}
\psi_{L_0,2} & \psi_{L_0,3} &\cdots &\psi_{L_0,L_0-1} &\psi_{L_0,L_0} &0 \\
\psi_{L_0,3} & \psi_{L_0,4} &\cdots  &\psi_{L_0,L_0} &0  &0  \\
\vdots & \vdots &\vdots &\vdots &\vdots &\vdots \\
\psi_{L_0,L_0} & 0  &\cdots & 0 &0 &0 \\
0 & 0  &\cdots & 0 &0  &0
\end{bmatrix}
+
\begin{bmatrix}
0 & 0  &\cdots & 0 &0  &0 \\
\psi_{L_0,1} & 0 &\cdots &0 &0 &0 \\
\vdots & \vdots &\vdots &\vdots &\vdots &\vdots \\
\psi_{L_0,L_0-2} & \psi_{L_0,L_0-3} &\cdots  &\psi_{L_0,1} &0  &0  \\
\psi_{L_0,L_0-1} & \psi_{L_0,L_0-2} &\cdots  &\psi_{L_0,2} &\psi_{L_0,1}  &0  \\
\end{bmatrix} .
\end{align*}
}
We thus obtain $\rho_{L_0} = - (I + \Phi_{L_0} )^{-1} \Psi_{L_0}$, $\gamma_0 = \sigma^2 / (1+\rho_{L_0}^\T \Psi_{L_0})$, and $\gamma_{\ell} = \gamma_0 \rho_{L_0,\ell} $ ($\ell = 1,\ldots,L_0$).
Furthermore, for each $\ell > L_0$, $\gamma_{\ell}$ equals to $-\sum_{k=1}^{L_0} \Psi_{L_0,k} \gamma_{\ell-k}$. 

\textbf{Case 2: }
The second case is AR$(L_0(N))$ with $L_0(N) = \lfloor N^{0.4} \rfloor$ and $ \Psi_{L_0(N)} = [0.7^k]_{k=1}^{L_0(N)}$, namely
$x_n+\psi_{L_0(N),1}x_{n-1}+ \cdots + \psi_{L_0(N),L_0(N)}x_{n-L_0(N)}=\epsilon_n, \, \epsilon_n \sim \mathcal{N}(0,1)$.
This is the case where the true order is large in terms of sample size, and thus it can be treated as the infinite dimensional model (see Remark~\ref{remark:ass}).
Note that all the roots of each characteristic polynomial have modulus $0.7$. 
For each sample size $N=100,500,1000,10000$, the true order that generated the autoregression is $6, 12, 15, 39$, respectively.

%

\textbf{Case 3: }
The third case is the first order moving average  process $x_n = \epsilon_n - 0.8 \epsilon_{n-1},\epsilon_{n}\sim \mathcal{N}(0,1) $. It is an autoregression with infinite order.
The exact mismatch error of an estimated filter $\Lambda_{L}$ could be calculated in the following way:
$\norm{\Lambda_{L} - \Psi_{\infty}}_{\Gamma_{\infty}}^2 
= E \bigl \{ x_{n+1} +  [x_n,\ldots,x_{n-L+1}] \Lambda_{L}  \bigr\}^2 \bigr \} - \sigma^2 
= 1.64 (1+\norm{\Lambda_{L}}_2^2 ) - 2 \cdot 0.8 \ ( \Lambda_{L,1} + \sum_{k=1}^{L-1} \Lambda_{L,k} \Lambda_{L,k+1} ) - 1
$, 
where we have used $E(x_n^2) = 1.64$, $E(x_n x_{n-1}) = -0.8$, and $E(x_n x_{n-k}) = 0$ for $k>1$. 

In summary, Table~\ref{table:result1} and \ref{table:result2} show that BC achieves the performance that we had expected: it is consistent when the model class is well-specified, and 
its predictive performance is always close to the optimum of AIC and BIC in both well-specified and mis-specified cases.
In practice when no prior knowledge about the model specification is available, the proposed method is more flexible and reliable than AIC and BIC in selecting the most appropriate dimension.


\begin{table}[h!]
    \caption{\label{table:result2} Mismatch errors (and their standard errors) of autoregressive models selected by BC, AIC, and BIC, along with the parametricness index, in three different cases (values except $\PI$ and its standard errors were rescaled by $10^{3}$)}
  \centering
\begin{tabular}{|c|cccc|cccc|}
\hline
\multirow{2}{*}{Case} & \multicolumn{4}{c|}{$N=$100}        & \multicolumn{4}{c|}{$N=$500}         \\
                      & BC     & AIC    & BIC    & $\PI$      & BC     & AIC    & BIC    & $\PI$       \\ \hline
\multirow{2}{*}{1}    & 19.7   & 28.6   & 16.6   & 0.96     & 2.9    & 5.7    & 2.4    & 0.97      \\
                      & (1.13) & (1.28) & (1.01) & (0.0061) & (0.18) & (0.26) & (0.13) & (0.0050)  \\ \hline
\multirow{2}{*}{2}    & 76.7   & 71.9   & 94.2   & 0.58     & 17.6   & 17.5   & 25.2   & 0.29      \\
                      & (1.24) & (1.08) & (1.33) & (0.016)  & (0.25) & (0.24) & (0.33) & (0.014)   \\ \hline
\multirow{2}{*}{3}    & 97.8   & 94.7   & 122.8  & 0.58     & 26.6   & 26.6   & 38.0   & 0.32       \\
					  & (1.28) & (1.12) & (1.55) & (0.016)  & (0.27) & (0.27) & (0.41) & (0.015)    \\ \hline
\multicolumn{1}{l}{}                        & \multicolumn{1}{l}{}   & \multicolumn{1}{l}{}   & \multicolumn{1}{l}{}   & \multicolumn{1}{l}{}          & \multicolumn{1}{l}{}   & \multicolumn{1}{l}{}  & \multicolumn{1}{l}{}   & \multicolumn{1}{l}{} 
%
\\ \hline
\multirow{2}{*}{Case} & \multicolumn{4}{c|}{$N=$1000}        & \multicolumn{4}{c|}{$N=$10000}     \\
                      & BC     & AIC    & BIC    & $\PI$      & BC     & AIC    & BIC    & $\PI$      \\ \hline
\multirow{2}{*}{1}    & 1.6    & 3.4    & 1.3     & 0.98     & 0.11    & 0.39    & 0.10    & 0.99     \\
                      & (0.11) & (0.15) & (0.065) & (0.0047) & (0.012) & (0.020) & 0.0049  & (0.0033) \\ \hline
\multirow{2}{*}{2}    & 9.9    & 9.9    & 14.6    & 0.18     & 1.4     & 1.4     & 2.1     & 0.11     \\
                      & (0.13) & (0.13) & (0.18)  & 0.012    & (0.019) & (0.019) & (0.025) & (0.0097) \\ \hline
\multirow{2}{*}{3}    & 14.6   & 14.6   & 22.1    & 0.21     & 2.02    & 2.02    & 3.19    & 0.032    \\
                      & (0.15) & (0.15) & (0.24)  & (0.013)  & (0.021) & (0.021) & (0.032) & (0.0056) \\ \hline
\end{tabular}
\end{table}

\subsection{Real data experiment: the El Nino data from 1935 to 2015}
\begin{figure}[h!]
\centering
  \includegraphics[width=1.0 \linewidth]{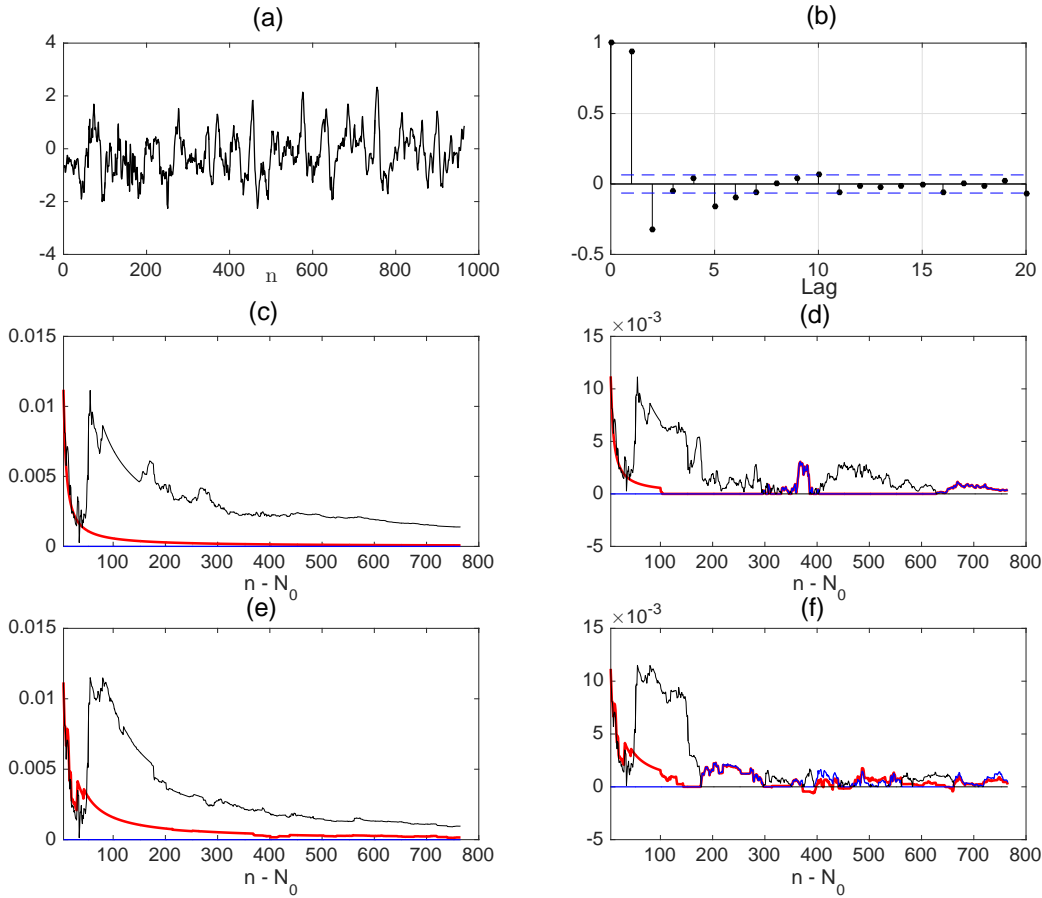}
  \vspace{-0.5in}
  \caption{(a) The monthly NINO3 index from January 1935 to May 2015; (b) the sample partial autocorrelations of the complete data with 95\% confidence bounds; (c) the normalized cumulated average prediction error at each time step (using all the current observations); (d) the normalized average prediction error over the recent window of size 100 (using all the current observations); (e) the normalized cumulated average prediction error (using the  recent $N_0$ observations); (f) the normalized average prediction error over the recent window of size 100 (using the  recent $N_0$ observations).
  In subfigures (c)-(f), BC, AIC, and BIC are respectively marked in red, blue, and black, and the curves have been normalized by subtracting the minimum of AIC curve and BIC curve.}
  \label{fig:EI_result}
  \vspace{-0.0in}
\end{figure}

As the largest climate pattern,  El Nino serves as the most dominant factor of oceanic influence on climate. The NINO3 index, defined as the area averaged sea surface temperature from $5^{\circ}$S-$5^{\circ}$N and $150^{\circ}$W-$90^{\circ}$W, is calculated from HadISST1 within the range of January 1935 to May 2015 \cite{rayner2003global}. The monthly data with overall 965 points is shown in Fig.~\ref{fig:EI_result}(a). The data seems to be highly dependent from its sample partial autocorrelations shown in Fig.~\ref{fig:EI_result}(b). 

To evaluate the predictive power of BC, AIC, and BIC, ideally we would apply  each estimated filter to independent and identically generated datasets as we have done in the synthetic data experiments. But it is not realistic to apply this cross-validation to a single real-world time series data. 
As an alternative, we adopt a prequential perspective \cite{dawid1984present,ing2005orderselection}, and evaluate the criteria in terms of the one-step prediction errors conditioning only on the past data at each time. 
Specifically, we start from an initial time step, say $N_0 = 200$, and obtain an estimated AR filter $\hat{\psi}_L(\mathcal{C})$ from the first $200$ points under each criterion $\mathcal{C}$. Upon the arrival of $(n=N_0+1)$th point, The one-step prediction error is revealed to be $\hat{e}_{n}(\mathcal{C}) = (x_{n} - [x_{n-1}, \ldots, x_{n-L}] \hat{\psi}_L )^2 $.  This procedure is repeated for $n=N_0+2,\ldots,N=965$, each time the AR filter being estimated from the observed $n-1$ data points and the tuning parameters being $\LM = \lfloor n^{1/3} \rfloor, \M = (\log n)^{0.9}$ (note that the $N$ in (\ref{eq:newer11}) was replaced by the available sample size $n$).
The cumulated average prediction error at each $n$ is computed to be $\overline{e}_{n}(\mathcal{C})=\sum_{t=N_0+1}^{n} \hat{e}_{t} (\mathcal{C}) / (n-N_0)$. 
To highlight the differences of $\overline{e}_{n}(\mathcal{C})$ for $\mathcal{C}=$ BC, AIC, BIC, we have plotted the normalized curve $\overline{e}_{n}(\mathcal{C}) - \overline{e}_{n}(opt)$ in Fig.~\ref{fig:EI_result}(c), where $\overline{e}_{n}(opt) = \min\{\overline{e}_{n}(\textrm{AIC}), \overline{e}_{n}(\textrm{BIC})\}$ for each $n=N_0+1,\ldots,N$.  
In order to show predictive power that may vary at different time epochs, We have also plotted in Fig.~\ref{fig:EI_result}(d) the (normalized) average prediction errors over only a sliding window of fixed size $100$, namely $\overline{e_0}_{n}(\mathcal{C}) = \sum_{t=s+1}^{n} \hat{e}_{t} (\mathcal{C}) / (n-s)$ where $s = \max \{N_0, n-100\}$.

In addition, in order to capture potential dynamics during different time epochs, we have also considered the estimation from a sliding window of fixed size $N_0$.
Specifically, we start from the same initial time step  $N_0 = 200$, and for each $n=N_0+1,\ldots,N$, the AR filters are estimated from only $n-N_0, \ldots, n-1$ with   $\LM = \lfloor N_0^{1/3} \rfloor, \M = (\log N_0)^{0.9}$ (note that the $N$ in (\ref{eq:newer11}) was replaced by the available sample size $N_0$). Similarly,  we computed the one-step prediction errors,  the normalized cumulated average prediction errors (plotted in Fig.~\ref{fig:EI_result}(e)), and the normalized windowed average prediction errors (plotted in Fig.~\ref{fig:EI_result}(f)). 
Fig.~\ref{fig:EI_result}(c)-(f) show that the performance of BC is close to AIC and outperforms BIC in general.

\subsection{Real data experiment: the English temperature data from 1659 to 2014}

\begin{figure}[h!]
\centering
  \includegraphics[width=1.0 \linewidth]{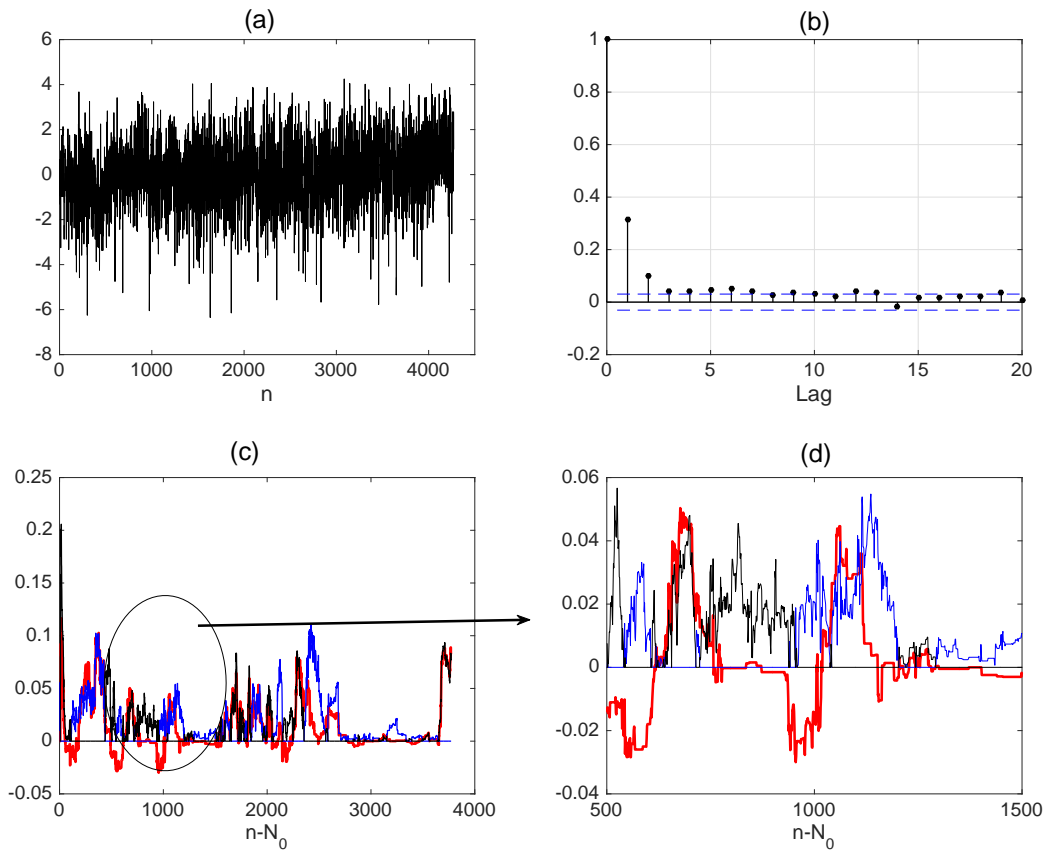}
  \vspace{-0.5in}
  \caption{(a) The de-seasoned data; (b) the sample partial autocorrelations of the complete data with 95\% confidence bounds; (c) the normalized cumulated average prediction error at each time step (using the recent $N_0$ observations); (d) the normalized average prediction error over the recent window of size 100 (using the  recent $N_0$ observations). In subfigures (c)-(d), BC, AIC, and BIC are respectively marked in red, blue, and black, and the curves have been normalized by subtracting the minimum of AIC curve and BIC curve.}
  \label{fig:CET_result}
  \vspace{-0.1in}
\end{figure}

In this experiment, we study the monthly English temperature data from 1659 to 2014 used by \citeasnoun{dieppois2013links}, which is perhaps the longest recorded environmental data in human history. 
We have pre-processed the raw data by subtracting each month by the average of that month over the 356 years. 
The de-seasoned data (with overall $N=4272$ points) is plotted in Fig.~\ref{fig:CET_result}(a). Its sample partial autocorrelations are shown in Fig.~\ref{fig:CET_result}(b).
In order to capture potential dynamics during such a long period, we adopt the prequential approach that was used to draw Fig.~\ref{fig:EI_result}(f), and omit the counterpart of Fig.~\ref{fig:EI_result}(c)(d)(e). 
Specifically, we started from $N_0=500$, and for each $n=N_0+1,\ldots,N$ the one-step ahead prediction was made by an AR filter produced from the recent window  of $N_0$ observations. The prediction errors $\hat{e}_{n} $ were averaged over a fixed window of size 100, namely $\overline{e_0}_{n}(\mathcal{C}) = \sum_{t=s+1}^{n} \hat{e}_{t} (\mathcal{C}) / (n-s)$ where $s = \max \{N_0, n-100\}$. 
We have plotted in Fig.~\ref{fig:CET_result}(c) the normalized average prediction errors, which is $\overline{e_0}_{n}(\mathcal{C}) - \overline{e_0}_{n}(opt)$ where $\overline{e_0}_{n}(opt) = \min\{\overline{e_0}_{n}(\textrm{AIC}), \overline{e_0}_{n}(\textrm{BIC})\}$ (similar as before).
We highlight the normalized average prediction errors  within the range $n=N_0+500,\ldots,N_0+1500$ in Fig.~\ref{fig:CET_result}(d). 
In this experiment, AIC is not constantly superior to BIC, and BC adaptively chooses to be close to the optimum of AIC and BIC. Furthermore, BC achieves the best predictive performance in some regions. 
The results show that BC is more flexible and reliable than AIC and BIC in practical applications. 
Note that we have adopted a specific choice of $\LM$ and $\M $ (see (\ref{eq:newer11})) throughout all the synthetic and real-world data experiments.
In practice, an analyst may achieve much better predictive performance of BC, by fine tuning $\LM$ and $\M $ for any particular real dataset. 

\section{Discussion}
\label{sec:Conclusion}


There have been many debates on which of AIC and BIC should be used \cite{burnham2004multimodel}. A practitioner who supports AIC may argue that all models are wrong, and thus it is safe to choose AIC that generally performs better in mis-specified situations. In contrast, a practitioner who supports BIC is usually in favor of the mathematically appealing ``consistency'' property and is quite confident that the candidate set of models contains the true (or practically a very good) model, or simply has a strong preference of parsimony in modeling.   
However, the debate is aroused due to the underlying assumption which tends to be overlooked: a practitioner should choose either AIC or BIC before even looking at the observed data---if some model specification test were done,  the practitioner might have changed his/her prejudice.
In a certain sense, the bent curve of bridge criterion, different from straight lines, was designed to mimic a sequence of model specification test which continuously check ``whether there exists a finite dimension $L_0$ underlying the observed data''.
For practical situations where there is no prior information, bridge criterion provides a practitioner with opportunities to change or reinforce his/her belief in the model specification. 

As a possible future work, it would be interesting to see in what extent the bridge criterion can be extended to other model selection problems, for instance 
the vector autoregressive model, autoregressive-moving-average model, and generalized linear  model. 

\section*{Acknowledgement}
This research was funded by the Defense Advanced Research Projects Agency (DARPA) under grant number W911NF-14-1-0508.
%


{
\small
\bibliographystyle{ECA_jasa} 
\bibliography{AR_order} 
}

\clearpage
\appendix

\end{document}